\documentclass{article}

\usepackage{def}
\usepackage{amssymb}
\usepackage{amsfonts}

\usepackage{amsmath}

\begin{document}

\title{Stable reduction of modular curves}

%\author[I.I.~Bouw]{Irene I.~ Bouw}
%\author[S.~Wewers]{Stefan Wewers}
%\address{%%
%Institut f\"ur Experimentelle Mathematik,\br
%Universit\"at GH Essen,\br
%Ellernstr. 29,\br
%45326 Essen, Germany}
%\email{bouw@exp-math.uni-essen.de}
%\address{%%
%Max--Plank--Institut f\"ur Mathematik,\br
%Vivatsgasse 7,\br
%53111 Bonn, Germany}
%\email{wewers@mpim-bonn.mpg.de}

\author{%
  Irene I.~ Bouw
\and
  Stefan Wewers
}

\date{}

\maketitle

\begin{abstract}
  We determine the stable reduction at $p$ of all three point covers
  of the projective line with Galois group ${\rm SL}_2(p)$. As a
  special case, we recover the results of Deligne and Rapoport on the
  reduction of the modular curves $X_0(p)$ and $X_1(p)$. Our method
  does not use the fact that modular curves are moduli
  spaces. Instead, we rely on results of Raynaud and the authors which
  describe the stable reduction of three point covers whose Galois
  group is strictly divisible by $p$.
\end{abstract}

\section*{Introduction}

Modular curves are quotients of the upper half plane by congruence
subgroups of ${\rm SL}_2(\ZZ)$. The most prominent examples are
$X_0(N)$, $X_1(N)$ and $X(N)$, for $N\geq 1$. Modular curves are also
moduli spaces for elliptic curves endowed with a level
structure. Therefore, they are defined over small number fields and
have a rich arithmetic structure.  Deligne and Rapoport \cite{DelRap}
determine the reduction behavior of $X_0(p)$ and $X_1(p)$ at the prime
$p$. Using this result, they prove a conjecture of Shimura saying
that the quotient of the Jacobian of $X_1(p)$ by the Jacobian of
$X_0(p)$ acquires good reduction over $\QQ(\zeta_p)$. Both the
reduction result and its corollary have been generalized by Katz and
Mazur \cite{KatzMazur} to arbitrary level $N$ and various level
structures.

The basic method employed in \cite{DelRap} and \cite{KatzMazur} is to
generalize the moduli problem defining the modular curve in question,
in such a way that it makes sense in arbitrary characteristic. 
By general results on representability of moduli problems one then
obtains a model of the modular curve over the ring of integers of a
subfield of $\QQ(\zeta_N)$. This model has bad reduction at all primes
dividing $N$, and it is far from being semistable, in general. In
spite of the very general results of \cite{KatzMazur}, it remains an
unsolved problem to describe the stable reduction of $X_0(N)$ at $p$
if $p^3$ divides $N$ (see \cite{Edixhoven90} for the case $p^2|| N$).

\vspace{1ex} In this note we suggest a different approach to study the
reduction of modular curves. Our starting point is the observation
that the $j$-invariant
\[
     X(N) \;\To\; X(1) \;\cong\; \PP^1
\]
presents the modular curve $X(N)$ as a Galois cover of the projective
line (with Galois group ${\rm PSL}_2(\ZZ/N)$) which is branched only
at the three rational points $0$, $1728$ and $\infty$. Let us call a
Galois covers of the projective line branched only at three points a
{\em three point cover}. In \cite{Raynaud98} Raynaud studies the
stable reduction of three point covers under the assumption that the
residue characteristic $p$ strictly divides the order of the Galois
group. His results have been sharpened in \cite{special} and
\cite{bad}. In the present paper we determine the stable reduction of
all three point covers with Galois group $\PSL_2(p)$, using the
results of \cite{bad}. As a corollary we obtain a new proof of the
results of Deligne and Rapoport on the reduction of $X_0(p)$ and
$X_1(p)$. Somewhat surprisingly, our proof does not use the modular
interpretation of these curves. On the other hand, the determination
of the stable reduction of some of the ${\rm PSL}_2(p)$-covers which
are not modular curves does use the fact that they are moduli spaces
of a certain kind. Here we build on the results of \cite{IreneAux}.

\vspace{1ex} What is so special about the group $\PSL_2(p)$? From our
point of view there are two main aspects. The first is {\em rigidity}.
Let $G$ be a finite group and ${\bf C}=(C_1,C_2,C_3)$ a triple of
conjugacy classes of $G$. Suppose that there exists a triple ${\bf
g}=(g_1,g_2,g_3)$ of generators of $G$ with $g_i\in C_i$ and
$g_1g_2g_3=1$. The triple ${\bf g}$ corresponds to a $G$-cover
$Y\to\PP^1$ with three branch points. If $G$ equals $\PSL_2(p)$, such
a triple ${\bf g}$ is unique, up to uniform conjugation in ${\rm
PGL}_2(p)$. Therefore there exists at most one three point cover with
a given branch cycle description ${\bf C}$, up to isomorphism.

We use rigidity as follows. By a result of \cite{bad} we can construct
three point $G$-covers $Y\to\PP^1$ with bad reduction at $p$ by
lifting a certain type of `stable $G$-cover' $\Yb\to\Xb$ from
characteristic $p$ to characteristic zero. Here $\Yb\to\Xb$ is a
finite map between semistable curves in characteristic $p$ together
with some extra structure $(g,\omega)$ which we call the {\em special
deformation datum}.  In the case of the modular curve $X(p)$ (where
$G={\rm PSL}_2(p)$) one can construct $\Yb\to\Xb$ very explicitly; the
special deformation datum $(g,\omega)$ corresponds to a solution of a
hypergeometric differential equation. (A similar phenomenon occurs in
Ihara's work on congruence relations \cite{Ihara74}.) The ${\rm
PSL}_2(p)$-cover $Y\to\PP^1$ resulting from the lifting process has
branch cycle description $(3A,2A,pA)$. By rigidity, $Y\to\PP^1$ is
isomorphic to $X(p)\to X(1)$. In particular, the stable reduction of
$X(p)$ is isomorphic to $\Yb$.

The other nice thing about ${\rm PSL}_2(p)$ is that $p$ strictly
divides its order. The results of \cite{Raynaud98} and \cite{bad}
which describe the stable reduction of a given $G$-cover require that
$p$ strictly divides the order of $G$, whereas the Sylow $p$-subgroup
of ${\rm PSL}_2(\ZZ/p^n)$ is rather big for $n>1$. This is the main
obstruction for extending the method of the present paper to modular
curves of higher $p$-power level. There are partial results of the
authors generalizing some of the results of \cite{Raynaud98} and
\cite{bad} to groups $G$ with a cyclic or an elementary abelian Sylow
$p$-subgroup (unpublished). It seems hopeless to obtain general
results beyond these cases. But maybe a combination of the methods
presented in the present paper with the modular approach might shed
some light on the stable reduction of modular curves of higher
$p$-power level.

The organization of this paper is as follows. In Section 1 we define
special deformation data and explain how to associate a special
deformation datum to a three point cover with bad reduction. We recall
a result from \cite{bad} which essentially says that we can reverse
this process. In Section 2 we introduce hypergeometric deformation
data. These are special deformation data satisfying an addition
condition (Definition \ref{hgdef}). We classify all hypergeometric
deformation data by showing that they correspond to the solution in
characteristic $p$ of some hypergeometric differential equation.  In
Sections 3 and 4 we use these results to give a new proof of the
stable reduction of the modular curves $X(p)$ and $X_0(p)$. In Section
5 we generalize these results to all three point covers with Galois
group $\SL_2(p)$.

\section{The special deformation datum}\label{sddsec}

Let $k$ be an algebraically closed field of characteristic $p>2$. Let
$H$ be a finite group of order prime to $p$. Fix a character
$\chi:H\to\FF_p^\times$.  Let $g:Z_k\to\PP^1_k$ be an $H$-Galois cover
and $\omega$ be a meromorphic differential on $Z_k$. We assume that
$\omega$ is logarithmic (i.e.\ can be written as $\omega={\rm d}u/u$)
and
\[
    \beta^\ast \omega=\chi(\beta)\cdot \omega, \qquad 
       \mbox{for all } \beta\in H.
\]
Let $\xi\in Z_k$ be a closed point and $\tau$ its image in
$\PP^1_k$. Denote by $H_\xi$ the stabilizer of $\xi$ in $H$. Define
\begin{equation}\label{definveq}m_\tau:=|H_\xi|, \qquad 
    h_\tau:=\Ord_\xi(\omega)+1, \qquad
     \sigma_\tau=h_\tau/m_\tau.
\end{equation} 
Since $\omega$ is logarithmic, we have $h_\tau\geq 0$.  We say that
$\tau$ is a {\em critical point} of the differential $\omega$ if
$(m_\xi, h_\xi)\neq (1,1)$. Let $(\tau_i)_{i\in\BB}$ be the critical
points of $\omega$, indexed by a finite set $\BB$. For $i\in\BB$, we
write $m_i, h_i, \sigma_i$ instead of $m_{\tau_i}, h_{\tau_i},
\sigma_{\tau_i}$.  For every $i$, choose a point $\xi_i\in Z_k$ above
$\tau_i$ and write $H({\xi_i})\subset H$ for its stabilizer.  Define
$\Bw:=\{\,i\in\BB \mid h_i=0 \,\}$.

\begin{defn}\label{sdddef}
  A {\em special deformation datum} of type $(H,\chi)$ is a pair $(g,
  \omega)$, where $g:Z_k\to \PP^1_k$ is an $H$-Galois cover and
  $\omega$ is a logarithmic differential on $Z_k$ such that the
  following holds.
\begin{itemize}
\item[(i)]  We have
\begin{equation}\label{actioneq}\beta^\ast \omega=\chi(\beta)\cdot \omega, \qquad \mbox{for all } \beta\in H.
\end{equation}
\item[(ii)] For every $i\in\BB-\Bw$, we have that $0<\sigma_i\leq 2$.
\item[(iii)] Define $\Bp=\{i\in\BB\, |\, 0<\sigma_i\leq 1\}$. Then
  $|\Bp\cup \Bw|=3$.
\end{itemize}
\end{defn}

Write $\BB_0=\Bw\cup \Bp$ and $\Bn=\BB-\BB_0$. See \cite{bad}
for more details and an explanation of the terminology. If
$\sigma_i\neq 1,2$ for all $i$, the above definition coincides with
\cite[Definition 2.7]{bad}.

\begin{defn}\label{taildef} Let $G$ be a finite group.
  A {\em $G$-tail cover} is a (not necessarily connected) $G$-Galois
  cover $f_k:Y_k\to\PP^1_k$ such that $f_k$ is wildly branched at
  $\infty$ of order $pn$ with $n$ prime to $p$ and tamely branched at
  no more than one other point.  We say that $f$ is a {\em primitive
  tail cover} if it is branched at two points. Otherwise, we call $f$
  a {\em new tail cover}.
\end{defn}

If the group $G$ is understood, we talk about tail covers instead of
$G$-tail covers.  To a $G$-tail cover $f_k:Y_k\to\PP^1_k$ we associate its
{\em ramification invariant} $\sigma(f)=h/n$, where $h$ is the
conductor of $f_k$ at $\infty$ and $n$ the order of the prime-to-$p$
ramification, as in Definition \ref{taildef}. The ramification
invariant is the jump in the filtration of higher ramification groups
in the upper numbering.

Three point covers with bad reduction give rise to
a special deformation datum and  a set of tail covers. Essentially,
the stable reduction of the cover is determined by these data.
Proposition \ref{liftprop} states that given a special
deformation datum and a set of tail covers satisfying some
compatibility conditions, there exists a three point cover in
characteristic zero which gives rise to the given datum.

To be more precise, let $R$ be a complete discrete valuation ring with
fraction field $K$ of characteristic zero and residue field an
algebraically closed field $k$ of characteristic $p$. Let $G$ be a
finite group and let $f:Y\to X=\PP^1_K$ be a $G$-Galois cover branched
at three points $x_1, x_2, x_3$. We assume that the points
$x_1,x_2,x_3$ specialize to pairwise distinct points on the special
fiber $\PP^1_k$. Denote the ramification points of $f$ by $y_1,\ldots,
y_s$. We consider the $y_i$ as markings on $Y$. After replacing $K$ by
a finite extension, there exists a unique extension $(Y_R; y_i)$ of
$(Y; y_i)$ to a stably marked curve over $R$. The action of $G$
extends to $Y_R$; write $X_R$ for the quotient of $Y_R$ by $G$. The
map $f_R:Y_R\to X_R$ is called the {\em stable model} of $f$; its
special fiber $\fb:\Yb\to\Xb$ is called the {\em stable reduction } of
$f$, \cite[Definition 1.1]{bad}.

We say that $f$ has {\em good reduction} if $\fb$ is separable. This
is equivalent to $\Xb$ being smooth. If $f$ does not have good
reduction, we say it has {\em bad reduction}.

\begin{figure}
\begin{center}
\unitlength3.5mm
%%% Local Variables: 
%%% mode: latex
%%% TeX-master: "mcav"
%%% End: 

\begin{picture}(30,16)

\put(6,4){\line(1,0){16}}
\put(6,10){\line(1,0){16}}
\put(6,14){\line(1,0){16}}

\put(10,3){\line(0,1){4}}
\put(12,3){\line(0,1){4}}
\put(16,3){\line(0,1){4}}
\put(19,3){\line(0,1){4}}

\put(10,9){\line(0,1){6}}
\put(12,9){\line(0,1){6}}
\put(16,9){\line(0,1){6}}
\put(19,9){\line(0,1){6}}

\put(8,4){\circle*{0.3}}
\put(10,6){\circle*{0.3}}
\put(12,6){\circle*{0.3}}

\put(17,5){\makebox(1,1){$\cdots$}}
\put(17,12){\makebox(1,0){$\cdots$}}

\put(7,11.9){\makebox(0,1){$\vdots$}}
\put(14,11.9){\makebox(0,1){$\vdots$}}
\put(21,11.9){\makebox(0,1){$\vdots$}}

\put(1.5,11.9){\makebox(1,1){$\Yb$}}
\put(1.5,3.5){\makebox(1,1){$\Xb$}}
\put(2,10){\vector(0,-1){3.5}}

\put(10,1.5){\makebox(2,1)
       {$\underbrace{\hspace{9mm}}_{\Bp}$}}
\put(16,1.5){\makebox(3,1)
       {$\underbrace{\hspace{14mm}}_{\Bn}$}}

\put(8.1,4.7){\makebox(0,0){$\scriptstyle \tau_1$}}
\put(9,5.9){\makebox(0,0){$\scriptstyle \bar{x}_2$}}
\put(11,5.9){\makebox(0,0){$\scriptstyle \bar{x}_3$}}

\put(23,4){\makebox(0,0){$\Xb_0$}}

\end{picture}
\caption{\label{stablepic} The stable reduction of $f:Y\to\PP^1$}
\end{center}
\end{figure}

Suppose that $f$ has bad reduction and that $p$ strictly divides the
order of $G$. Then $\Xb$ is a `comb', see \cite[Theorem 2.6]{bad} and
Figure \ref{stablepic}.  The central component $\Xb_0\subset\Xb$ is
canonically isomorphic to $\PP^1_k$.  Choose a component $\Yb_0$ of
$\Yb$ above $\Xb_0$. Let $G_0\subset G$ denote the decomposition group
of the component $\Yb_0$ and $I_0\lhd G_0$ the inertia group. Then
$I_0\cong\ZZ/p$. Therefore, the restriction of $\fb$ to $\Yb_0$
factors as $\Yb_0\to\Zb_0\to\Xb_0$, with $g:\Zb_0\to\Xb_0$ a separable
Galois cover of prime-to-$p$ order and $\Yb_0\to\Zb_0$ purely
inseparable of degree $p$. The map $\Yb_0\to\Zb_0$ is generically
endowed with the structure of a $\bmu_p$-torsor. This structure is
encoded in a logarithmic differential $\omega$. Define
$H:=\Gal(\Zb_0/\Xb_0)$; then $G_0$ is an extension of $H$ by
$I_0$. The action of $H$ on $I_0$ by conjugation gives rise to a
character $\chi:H\to\FF_p^\times$. It is proved in \cite[Proposition
2.15]{bad} that $(g,\omega)$ is a special deformation datum of type
$(H,\chi)$.

The irreducible components of $\Xb$ different from the central
component $\Xb_0$ are called the {\em tails} of $\Xb$.  The tails are
parameterized by $i\in\BB-\Bw$, \cite[Proposition1.7]{bad}. For
$i\in\BB-\Bw$, write $\Xb_i$ for the corresponding tail of $\Xb$ and
write $\fb_i:\Yb_i\to\Xb_i$ for the restriction of $\fb$ to $\Xb_i$.
The map $\fb_i$ is a primitive tail cover if $i\in\Bp$ and a new tail
cover if $i\in\Bn$, \cite[Lemma 1.4]{bad}. We say that $\Xb_i$
is a {\em primitive} (resp.\ {\em new}) tail of $\Xb$ if $i\in\Bp$
(resp.\ $i\in\Bn$). 

The three branch points specialize as follows. If $p$ divides the
ramification index of $x_j$, then $x_j$ specializes to
$\tau_i\in\Xb_0$ for a unique $i\in\Bw$. Otherwise, $x_j$ specializes
to one of the primitive tails $\Xb_i$, for  $i\in\Bp$. From now on, we will
identify the set $\BB_0=\Bp\cup\Bw$ with $\{1,2,3\}$ such that
$i\in\Bw$ if and only if $p$ divides the order of the ramification
index at $x_i$. For $i\in\Bp$ we denote by $\bar{x}_i\in\Xb_i$ the
specialization of $x_i$. 

The special deformation datum $(g, \omega)$ and the tail covers
$\fb_i$ satisfy the following compatibility condition. For every
$i\in\BB-\Bw$, there exists a point $\xi_i\in\Zb_0$ above $\tau_i$ and
a point $\eta_i\in\Yb_i$ which is wildly
ramified in $\fb_i$ such that
\begin{equation}\label{compeq}
  I(\eta_i)\subset G_0,\qquad I(\eta_i)\cap H=H({\xi_i}).
\end{equation}
Here $I(\eta_i)\subset G$ (resp.\ $H(\xi_i)\subset H$) denotes the
stabilizer of $\eta_i$ (resp.\ of $\xi_i$). Furthermore, the
ramification invariant $\sigma(\fb_i)$ of the tail cover $\fb_i$ is
equal to $\sigma_i$ as in Definition \ref{sdddef}.

\begin{prop}\label{liftprop}
  Let $G$ be a finite group. Fix a subgroup $G_0\simeq
  \ZZ/p\rtimes_\chi H$ of $G$, where $H$ has order prime-to-$p$ and
  acts on $\ZZ/p$ via a character $\chi:H\to\FF_p^\times$.  Let $(g,
  \omega)$ be a special deformation datum of type $(H,\chi)$. Let
  $\BB$, $\tau_i$ and $\sigma_i=h_i/m_i$ be as in Definition
  \ref{sdddef}.
  
  For $i\in\Bp$ (resp.\ $i\in\Bn$), suppose we are given a primitive
  (resp.\ new) $G$-tail cover $\fb_i:\Yb_i\to \Xb_i=\PP^1_k$ with
  ramification invariant $\sigma_i$. Suppose furthermore that there
  exists a point $\xi_i\in Z_k$ above $\tau_i\in\PP^1_k$ and a point
  $\eta_i\in \Yb_i$ which is wildly ramified in $\fb_i$ such that
  \eqref{compeq} holds.
  
  Let $K_0$ be the fraction field of $W(k)$ (the Witt vectors over
  $k$) and let $K/K_0$ be the (unique) tame extension of degree
  \begin{equation} \label{Neq}
      N \;:=\; (p-1)\cdot \mathop{\rm lcm}_{i\not\in\Bw}(h_i).
  \end{equation}
  There exists a $G$-Galois cover $f:Y\to X=\PP^1_K$ over $K$,
  branched at $0$, $1$ and $\infty$, which gives rise to the special
  deformation datum $(g, \omega)$ and the tail covers $\fb_i$.
\end{prop}

\begin{Proof}
  The datum $(g,\omega, \fb_i)$ together with the choice of the points
  $\xi_i$ and $\eta_i$ is called a {\em special $G$-deformation datum}
  in \cite[Definition 2.13]{bad}. Hence the proposition follows from
  \cite[Corollary 4.6]{bad}.
\end{Proof}

In this paper, the group $G$ is either $\SL_2(p)$ or
  $\PSL_2(p)$. Therefore, the group $H$ is cyclic and $H(\xi_i)$ only
  depends on $i$ and not on the choice of $\xi_i$.

The $G$-cover $f:Y\to\PP^1_K$ in Proposition \ref{liftprop}
  which gives rise to the datum $(g,\omega,\fb_i)$ is not unique.
In  \cite[Theorem 4.5]{bad} one finds an explicit parameterization of the
  the set of isomorphism classes of all such $G$-covers.

\section{Hypergeometric deformation data }\label{hgsec}

In this section we classify a certain type of special deformation data
which we call {\em hypergeometric}. A hypergeometric (special)
deformation datum corresponds to a polynomial solution of a
hypergeometric differential equation. For this reason it is possible
to find all such deformation data. The definition of hypergeometricity
looks very restrictive, but we show in Section \ref{legosec} that the
special deformation datum corresponding to a three point cover with
Galois group $\PSL_2(p)$ or $\SL_2(p)$ is hypergeometric. In
particular, this holds for the $\PSL_2(p)$-cover $X(p)\to
X(1)\simeq\PP^1_j$.

\begin{defn}\label{hgdef}
  A {\em hypergeometric deformation datum} is a special deformation
  datum $(g, \omega)$ with $\sigma_i=(p+1)/(p-1)$ for all
  $i\in\Bn$.
\end{defn}

Let $(g, \omega)$ be a hypergeometric deformation datum.  Recall that
$\BB_0=\{1,2,3\}$, by assumption. In this section we assume $x_1=0,
x_2=1, x_3=\infty$. In particular, $\tau_i\neq \infty$ for all
$i\in\Bn$. For $i\in\BB$, we define integers $0\leq a_i< p-1$ such
that
\[a_i\equiv (p-1)\sigma_i\mod{p-1}.\]
In particular, $a_i=2$ for $i\in \Bn$.

 To $g:Z_k\to\PP^1_k$ we associate the (not necessarily connected)
$\FF_p^\times$-Galois cover
\[
        g':Z_k':=\Ind_{\Im(\chi)}^{\FF_p^\times} (Z_k/{\Ker{\chi}}) 
          \to \PP^1_k.
\]
It follows easily from Definition \ref{hgdef} and the expansion of
$\omega$ into local coordinates that $Z_k'$ is the complete
nonsingular curve defined by
\begin{equation}\label{typeeq}z^{p-1}= x^{a_1}(x-1)^{a_2}\prod_{i\in\Bn} (x-\tau_i)^2,\end{equation}
where $\beta^\ast z= \chi(\beta)\cdot z$. We call $(a_1, a_2,
a_3)$ the {\em signature} of the hypergeometric deformation datum
$(g, \omega)$.

We may identify the differential $\omega$ on $Z_k$ with a differential
on $Z_k'$ which we also denote by $\omega$. Definition \ref{sdddef}(i)
and (\ref{definveq}) imply that
\begin{equation}\label{omegaeq}
       \omega=\epsilon\frac{z\, {\rm d} x}{x(x-1)},
\end{equation}
for some $\epsilon\in k^\times$.  The condition that $\omega$ is
logarithmic imposes a strong condition on $u:=\prod_{i\in
  \Bn}(x-\tau_i)$. Proposition \ref{hgprop} below states that the
differential (\ref{omegaeq}) is logarithmic if and only if $u$ is the
solution to a certain hypergeometric differential equation.

It follows from (\ref{definveq}) and the Riemann--Roch
Theorem applied to $\omega$ that 
\begin{equation}\label{vceq}
|\Bn|=(p-1-a_1-a_2-a_3)/2.
\end{equation}
In particular, $a_1+a_2+a_3$ is an even integer between 0 and
$p-1$.  Note that if $a_1+a_2+a_3=p-1$ then
$\deg(u)=|\Bn|=0$. Therefore, we exclude this case in Proposition
\ref{hgprop}.(i). However, Proposition \ref{hgprop}.(ii) applies also
if $a_1+a_2+a_3=p-1$.

\begin{prop}\label{hgprop}
\begin{itemize}
\item[(i)] Let $(g,\omega)$ be a hypergeometric deformation datum. Let
  $(a_1, a_2, a_3)$ be its signature, and suppose that
  $a_1+a_2+a_3<p-1$.  Then $u(x):=\prod_{i\in\Bn} (x-\tau_i)$
is a solution to the hypergeometric differential equation
\[x(x-1) u''+[(A+B+1)x-C]u'+ABu=0,\]
where $A=(1+a_1+a_2+a_3)/2$, $B=(1+a_1+a_2-a_3)/2$ and
$C=1+a_1$.  The degree of $u$ in $x$ is $d:=(p-1-a_1-a_2-a_3)/2$.
\item[(ii)] Let $0\leq a_1, a_2, a_3<p-1$ be integers with
$a_1+a_2+a_3\leq p-1$. We assume that $a_1+a_2+a_3$ is
even. Let $H=\ZZ/m$, where $m=(p-1)/2$ if $a_1, a_2, a_3$ are all
even and $m=p-1$ otherwise. Choose an injective character
$\chi:H\to\FF_p^\times$.  There exists a hypergeometric deformation
datum $(g, \omega)$ of signature $(a_1, a_2, a_3)$ and type
$(H,\chi)$. The deformation datum $(g, \omega)$ is uniquely determined
by $(a_1, a_2, a_3)$, up to multiplying the differential $\omega$
by an element of $\FF_p^\times$.
\end{itemize}
\end{prop}

\begin{Proof}
  The following proof is inspired by \cite[Lemma 3]{Beukers02}. A
  similar argument can be found in \cite[Theorem 5]{Ihara74}. Write
\[
   Q=x^{1+a_1}(x-1)^{1+a_2}, \qquad 
           F=\epsilon\frac{z}{x(x-1)}=\frac{\epsilon z^p}{Qu^2}.
\]
It is well known that the fact that $\omega=F\, {\rm d}x$ is
logarithmic is equivalent to $D^{p-1} F=-F^p,$ where
$D:={\rm d}/{\rm d}x$.  Since $D^{p-1} F=\epsilon
z^pD^{p-1}[1/(Qu^2)]$, we find
\begin{equation}\label{Cartiereq}
  D^{p-1}\frac{1}{Qu^2}= -\epsilon^{p-1}\frac{1}{x^p(x-1)^p}.
\end{equation}

Choose $i\in\Bn$ and write
\[
   \frac{1}{Qu^2} \;=\; \sum_{n\geq -2} a_n (x-\tau_i)^n.
\]
Then
\[  
   D^{p-1}\frac{1}{Qu^2} \;=\; -\left[\frac{a_{-1}}{(x-\tau_i)^p}+
        a_{p-1}+\cdots\right] \;=\; -\epsilon^{p-1}\frac{1}{x^p(x-1)^p}.
\]
We conclude that $a_{-1}=0$. 

Write
\begin{eqnarray*}
   Qu^2 &=&  [Q(\tau_i)+Q'(\tau_i)(x-\tau_i)\cdots]
   [u'(\tau_i)(x-\tau_i)+\frac{1}{2}u''(\tau_i)(x-\tau_i)^2+\cdots]^2\\
        &=&  (x-\tau_i)^2\,[\,Q(\tau_i)u'(\tau_i)+(\,Q'(\tau_i)u'(\tau_i) 
               + Q(\tau_i)u''(\tau_i)\,)\,(x-\tau_i)+\cdots].
\end{eqnarray*}
We see that
$-a_{-1}=u'(\tau_i)[\,Q'(\tau_i)u'(\tau_i)+Q(\tau_i)u''(\tau_i)\,]$. Since
$u'(\tau_i)\neq 0$, we have that
$Q'(\tau_i)u'(\tau_i)+Q(\tau_i)u''(\tau_i)=0$ for all $i\in \Bn$.

Define $G=Q'u'+Qu''$. This is a polynomial of degree less than or
equal to $e:=\deg(Q)+\deg(u)-2=|\Bn|+a_1+a_2$. The coefficient of
$x^e$ in $G$ is $g_e=\deg(u)(\deg(Q)+\deg(u)-1)$. Since $\deg(Q)\geq
2$ and $\deg(u)\neq 0$, the coefficient $g_e$ is nonzero.  The
polynomial $G$ is obviously divisible by $x^{a_1}(x-1)^{a_2}$, and
$G(\tau_i)=0$, for all $i\in\Bn$. Therefore $G=g_e x^{a_1}(x-1)^{a_2}
u$. Dividing by $x^{a_1}(x-1)^{a_2}$, we find that $u$ is a solution to
the hypergeometric differential equation
\begin{equation}\label{hgeq}
  P_0 u''+P_1 u'+P_2 u=0,\end{equation} where $P_0=x(x-1)$,
$P_1=P_0Q'/Q=(2+a_1+a_2)x-(1+a_1)$, and
$P_2=-g_e=((1+a_1+a_2)^2-a_3^2)/4$. Part (i) of the proposition
follows. The proof also shows that 
\begin{equation} \label{epsiloneq}
  \epsilon^{p-1} \;=\; (-1)^{a_1} 
      \begin{pmatrix} p-1-a_2 \\ a_1 \end{pmatrix}.
\end{equation}

To prove (ii), let $0\leq a_1, a_2, a_3<p-1$ be integers such that
$a_1+a_2+a_3$ is even and less than or equal to $ p-1$. Define $d=(p-1-a_1-a_2-a_3)/2$. We
claim that the differential equation (\ref{hgeq}) has a unique
polynomial solution $u=\sum_i u_i x^i\in k[x]$ of degree $d$ with
$u_d=1$.

Suppose that $u=\sum_{i\in\ZZ} u_i x^i$ is a solution of (\ref{hgeq}).
Put $A_i=(i+1)(i+a_1+1)$ and
$B_i=(i+(1+a_1+a_2+a_3)/2)(i+(1+a_1+a_2-a_3)/2)$.  The
differential equation gives a recursion
\begin{equation}\label{recursioneq}A_i u_{i+1}= B_iu_i\end{equation}
for the coefficients of $u$. One checks that $A_i$ and $B_i$ are
nonzero for all $0\leq i < d$ and that $A_{-1}=B_d=0$ and $A_d\cdot
B_{-1}\neq 0$. Therefore there is a unique polynomial solution $u$ of
degree $d$ with $u_d=1$.  This implies the claim.

Let $m$ be as in the statement of the proposition. Let $Z_k'$ be the
complete nonsingular curve defined by
\[z^{p-1}=u^2x^{a_1}(x-1)^{a_2}\]
and let $Z_k$ be a connected component of $Z'_k$.  Write
$g:Z_k\to\PP^1_k$ for the $m$-cyclic cover defined by $(x,z)\mapsto
x$. We may identify $H$ with $\Gal(Z_k/\PP^1_k)$ such that $H$
acts on $z$ via some chosen character $\chi:H\to\FF_p^\times$. Let
\[
   \omega=\frac{\epsilon\, z\,{\rm d} x}{x(x-1)}.
\]

 Since $u$ is the solution to a hypergeometric differential equation,
it has at most simple zeros outside 0 and 1. We know that
$u(0)=u_0=\prod_{i=0}^{d-1} (B_i/A_i)\neq 0$. Therefore $u$ does not
have a zero at $x=0$. We claim that $u$ does not have a zero at 1,
also. To see this, let $t=1-x$ be a new coordinate on $\PP^1_k$ and
write $\tilde{u}(t)=u(1-t).$ The coefficients of $\tilde{u}$ satisfy a
recursion similar to  \eqref{recursioneq}, but with $a_1$ and $a_2$
interchanged. The same argument as above applied to $\tilde{u}$ shows
now that $\tilde{u}(0)=u(1)\neq 0$.  We conclude that $u$ has exactly
$d$ zeros different from 0 and 1.  Analogous to the proof of (i), one
shows that for $\epsilon$ as in \eqref{epsiloneq} the differential
$\omega$ is logarithmic. It is clear that $(g, \omega)$ is a special
deformation datum.  The uniqueness statement is obvious.
\end{Proof}

\section{The reduction of $X(p)$} \label{modularsec}

In the rest of the paper we suppose that $p\geq 5$. Choose a primitive
$(p-1)$th root of unity $\zeta\in \FF_p$ and a primitive $(p+1)$th
root of unity $\tilde{\zeta}\in \FF_{p^2}$. Define
\[
    C(l) \;=\; \{ A\in \SL_2(p)\,|\, \tr(A)=\zeta^l+\zeta^{-l}\},\qquad
        \mbox{for } 0<l\leq (p-1)/2
\]
and
\[
  \tilde{C}(l) \;=\; \{ A\in \SL_2(p)\,|\,\tr(A)=\tilde{\zeta}^l+
       \tilde{\zeta}^{-l}\},\qquad \mbox{for } 0<l< (p+1)/2.
\]
These are the conjugacy classes of $\SL_2(p)$ of nontrivial elements
of order prime to $p$. We write $pA$ and $pB$ for the two conjugacy
classes of elements of order $p$. Note that the elements of order $p$
have trace 2.

Let ${\bf C}=(C_1, C_2, C_3)$ be a triple of conjugacy classes of
$\SL_2(p)$. In the rest of this section, we suppose that the $C_i$ do
not contain $\pm I$.  We write $\Niin_3({\bf C})$ for the set of
isomorphism classes of $\SL_2(p)$-covers $Y\to \PP^1$ branched at $0,
1, \infty$ with class vector ${\bf C}$. This means that the canonical
generator of some point of $Y$ above $x_i\in \{0, 1, \infty\}$ with
respect to the chosen roots of unity is contained in the conjugacy
class $C_i$. In this paper, we call two $G$-covers $f_i:Y_i\to X$
isomorphic if there exists a $G$-equivariant automorphism $\phi:Y_1\to
Y_2$ such that $f_1=f_2\circ \phi$.

\begin{prop}[Linear rigidity]\label{linrigprop} Let ${\bf C}=(C_1, C_2, C_3)$ be a triple of conjugacy
classes of $\SL_2(p)$.
\begin{itemize}
\item[(i)] Suppose that the elements of $C_i$ have prime-to-$p$
  order. Then $ \# \Niin_3({\bf C})\in\{0, 2\}.$
\item[(ii)] 
  Otherwise, $ \# \Niin_3({\bf C})\in\{0, 1\}.$
\end{itemize}
\end{prop}

\begin{Proof}
  See \cite{StrambachVoelklein99}. The difference between the two
  cases comes from the fact that the outer automorphism group of $G$
  has order two and interchanges the two conjugacy classes of order
  $p$ but fixes all other conjugacy classes.
\end{Proof}

In the following we denote by $2A$ (resp.\ $3A$) the unique conjugacy
classes in ${\rm PSL}_2(p)$ of elements of order $2$ (resp.\ $3$).
Note that elements of $2A$ (resp.\ elements of $3A$) lift to elements
of order $4$ (resp.\ of order $3$ or $6$) in ${\rm SL}_2(p)$.
Furthermore, we denote by $pA$ and $pB$ the images in ${\rm PSL}_2(p)$
of the conjugacy classes of ${\rm SL}_2(p)$ with the same name.
 
Let $X(N)$ be the modular curve parameterizing (generalized)
elliptic curves with full level $N$-structure \cite{DelRap}.  Consider
the cover $X(2p)\to X(2)\simeq \PP^1_\lambda$.  This is a Galois cover
with Galois group $\PSL_2(p)$ which is branched at $0,1,\infty$ of
order $p$ and unbranched elsewhere.  One easily checks as in
\cite[Lemma 3.27]{Voelklein} that there are no $\PSL_2(p)$-covers with
class vector $(pA, pA, pB)$.  After renaming the conjugacy classes,
the class vector of $X(2p)\to X(2)$ is therefore ${\bf C}=(pA, pA,
pA)$.

Proposition \ref{linrigprop} implies that the triple ${\bf C}=(pA, pA,
pA)$ is rigid, i.e.\ there is a unique $\PSL_2(p)$-cover of
$\PP^1$ whose class vector is ${\bf C}$, up to automorphism. Since
$\PSL_2(p)$ has trivial center, it follows that the ${\rm
  PSL}_2(p)$-cover $X(2p)\to X(2)$ has a unique model $X(2p)_K\to
X(2)_K=\PP^1_K$ over any field $K$ of characteristic zero containing
$\sqrt{p^*}$, \cite[Chapter 7]{SerreTopics}.  Let $K_0$ be the
fraction field of the ring of Witt vectors over $k=\bar{\FF}_p$, and
let $K/K_0$ be a sufficiently large finite extension.  We define the
stable reduction of $X(2p)\to X(2)$ at $p$ as the stable reduction of
$X(2p)_K\to\PP^1_K$.

The cover $X(2p)\to X(2)$ may be considered as a variant with ordered
branch points of the cover $X(p)\to X(1)\simeq \PP^1_j$. We find it
easier to first compute the stable reduction of the former and then
deduce the stable reduction of the latter.

\begin{prop}\label{redmodprop}
  Let $\Yb\to\Xb$ be the stable reduction of $X(2p)\to X(2)$ at $p$.
  Then $\Xb$ has $(p-1)/2$ new tails with ramification invariant
  $\sigma=(p+1)/(p-1)$. The points $\tau_i\in\PP^1_k$, $i\in\Bn$ (where
  the new tails occur) are precisely the zeros of the polynomial
  \begin{equation} \label{ueq}
      u  \;=\; \sum_{j=0}^{(p-1)/2}\binom{(p-1)/2}{j}^2 x^j.
  \end{equation}
  Moreover, the stable reduction occurs over the (unique) tame
  extension $K/K_0$ of degree $(p^2-1)/2$. 
\end{prop}

The polynomial $u$ is known as the {\em Hasse invariant}.  It has the
property that the elliptic curve given by $y^2=x(x-1)(x-\lambda)$ is
supersingular in characteristic $p$ if and only if $\lambda$ is a zero
of $u$.  It is a solution to  the hypergeometric
differential equation with parameters $(1/2, 1/2, 1)$.  The relation
between this differential equation and the reduction of elliptic curves
is well understood, see for example \cite{Katz84}.

\begin{Proof}
  Let $a_1=a_2=a_3=0$ and $m=(p-1)/2$. Choose an injective character
  $\chi:\ZZ/m\to\FF_p^\times$. Let $(g, \omega)$ be the special
  deformation datum of signature $(0,0,0)$ constructed in Proposition
  \ref{hgprop}.(ii). Recall that we associate to the triple $(a_1,
  a_2, a_3)=(0,0,0)$ the hypergeometric differential equation with
  parameters $A=B=1/2$ and $C=1$. Using the recursion
  (\ref{recursioneq}) one checks that the polynomial $u$ defined by
  \eqref{ueq} is a solution in characteristic $p$ to this differential
  equation.  This differential equation has a unique monic solution of
  degree $(p-1)/2$.  Therefore, $u$ is the same as in Proposition
  \ref{hgprop} (ii).  
  
  Let $Y'_k\to\PP^1_k$ be a connected and new $\PSL_2(p)$-tail cover
  with ramification invariant $\sigma=(p+1)/(p-1)$. The existence of
  such a tail cover is proved in Lemma \ref{AL} below. The inertia
  group of a point above $\infty$ has order $p(p-1)/2$.  For every
  $i\in\Bn$, we let $\fb_i:\Yb_i\to\Xb_i$ be a copy of
  $Y_k'\to\PP^1_k$.  Let $\eta_i\in\Yb_i$ be the point whose inertia
  group consists of the upper triangular matrices in $\PSL_2(p)$. Then
  the condition of Proposition \ref{liftprop} is satisfied.
  Therefore, there exists a $\PSL_2(p)$-cover $f_K:Y_K\to\PP^1_K$
  branched at three points of order $p$ whose stable reduction gives
  rise to the special deformation datum $(g, \omega)$ and the tail
  covers $\fb_i$. This cover is defined and has stable reduction over
  the tame extension $K/K_0$ of degree $(p^2-1)/2$.  It follows from
  rigidity that $f_K$ is isomorphic to $X(2p)_K\to X(2)_K$. This
  proves the proposition.
\end{Proof}

\begin{lem}\label{AL} 
  There exist a connected  new $\SL_2(p)$-tail cover $h: Y_k\to
  X_k=\PP^1_k$, with inertia group of order $p(p-1)$ and ramification
  invariant $\sigma=(p+1)/(p-1)$.
\end{lem}

\begin{Proof}
Let $Y_k$ be the complete nonsingular curve defined by
\[
       y^{p+1} \;=\; x^p-x.
\]
The group $\SL_2(p)$ acts on $Y_k$ as follows. For
\begin{gather*}A=\begin{pmatrix}
a&b\\ c&d
\end{pmatrix}\in \SL_2(p), \quad \mbox{ define } Ax=\frac{ax+b}{cx+d}, \quad Ay=\frac{y}{cx+d}.\end{gather*}
It is easy to check that the quotient map $h:Y_k\to X_k\simeq
\PP^1_k$ is branched at exactly one point which we may assume to be
$\infty$. The inertia group of some point above $\infty$ has order
$p(p-1)$ and $\sigma=(p+1)/(p-1)$.
\end{Proof}

Let $X(p)\to X(1)\simeq \PP^1_j$ be the projection of the modular
curve $X(p)$ to the $j$-line. This is a $\PSL_2(p)$-cover with branch cycle description
$(3A,2A,pA)$ , which we may assume to be  branched at
$x_1=0, x_2=1728, x_3=\infty$. By Proposition \ref{linrigprop}, it has a unique model $X(p)_K\to\PP^1_K$
over any field $K$ of characteristic $0$ containing $\sqrt{p^*}$. As
before, we take $K/K_0$ a sufficiently large finite extension and
define the stable reduction of $X(p)\to X(1)$ at $p$ as the stable
reduction of the model $X(p)_K\to\PP^1_K$.

\begin{cor} \label{xpcor}
  Define $\alpha=\lfloor p/12\rfloor$.  Let $\Yb'\to\Xb'$ be the stable
  reduction of $X(p)\to X(1)$. Then $\Xb'$ has two primitive tails and
  $\alpha$ new tails. For $i=1,2$, denote by $\Xb_i$ the primitive
  tail to which $x_i$ specializes.  Then
\[
   \sigma_i \;=\; \left\{\;\;\; \renewcommand{\arraystretch}{2}
  \begin{array}{ll} 
  \displaystyle{\frac{p+1}{p-1}}& \mbox{\rm if\ }\; i\in \Bn,\\
  \displaystyle{\frac{\lceil(p-1)/3\rceil}{p-1}}&\mbox{\rm if\ }\; i=1,\\
  \displaystyle{\frac{2\lceil(p-1)/4\rceil}{p-1}}&\mbox{\rm if\ }\; i=2.
  \end{array}\right.
\]
Moreover, the stable reduction occurs over the tame extension $K/K_0$
of degree $(p^2-1)/2$. 
\end{cor}

\begin{Proof}
Consider the commutative diagram
\[\renewcommand{\arraystretch}{1.5}
\begin{array}{ccc}
  X(2p)_K &\To& X(p)_K\\
  \big\downarrow && \big\downarrow\\
  X(2)_K &\To& X(1)_K.
\end{array}
\]
The horizontal arrows are Galois covers whose Galois group is the
symmetric group $S_3$ on three letters. Let $Y_R$ be the stable model
of $X(2p)_K$. The action of $S_3$ extends to $Y_R$ and
$Y_R/S_3:={Y}^\ast_R$ is a semistable curve,
\cite[Appendice]{Raynaudfest}. 

Since $S_3$ acts on $Y_R$ it follows that the set of new tails of
$\Xb_0$ is stable under the action of $S_3$.  We claim that the action
of $S_3$ on the set of new tails of $\Xb$ has $\alpha$ orbits of
length six, one orbit of length three if $p\equiv 3\bmod 4$ and one
orbit of length two if $p\equiv 2\bmod 3$. Since $\Xb_0$ has genus
zero, an element of order three in $S_3$ has two fixed points in
$\Xb_0-\{0,1,\infty\}$. An element of order two in $S_3$ fixes one of
the points $0,1,\infty$ and has therefore one fixed point in
$\Xb_0-\{0,1,\infty\}$. The claim now follows by distinguishing the
four possibilities for $p\bmod{12}$.

The curve $X^\ast_R:=Y^\ast_R /G$ is almost equal to $X'_R$. The
special fiber $\Xb^\ast$ has a tail at zero (resp.\ 1728) if and only
if $p\equiv -1\bmod{3}$ (resp.\ $p\equiv -1\bmod{4}$), i.e.\ exactly
when zero (resp.\ 1728) is supersingular.  To obtain $X_R'$ one has to
perform a suitable blow up centered at the point $j=0$ (resp.\
$j=1728$) on $\Xb^\ast$ if $p\equiv 1\bmod 3$ (resp.\ $p\equiv 1\bmod
4$.)  The corollary follows, by analyzing the ramification of
$\Yb_0^\ast\to\Xb_0^\ast$. For future reference, we note that a
connected component of $\Yb_1'$ has decomposition group $\PSL_2(p)$ if
$p\equiv -1\bmod{3}$ and the nonabelian group of order $3p$ if
$p\equiv 1\bmod{3}$. A connected component of $\Yb_2'$ has
decomposition group $\PSL_2(p)$ if $p\equiv -1\bmod{4}$ and the
dihedral group of order $2p$ if $p\equiv 1\bmod{3}$.
\end{Proof}

The hypergeometric differential equation corresponding to $X(p)\to
X(1)$ is
\[
 x(x-1728)u''+[(2+a_1+a_2)x-1728(1+a_1)]u'+(1+a_1+a_2)^2u/4,
\] where
$a_1=\lceil(p-1)/3\rceil$ and $a_2=2\lceil (p-1)/4\rceil$. The
polynomial
\[
   u=(1728)^\alpha\sum_{n=0}^\alpha
   \binom{a_1+\alpha}{\alpha-n}\binom{\alpha}{n}\left(\frac{x}{1728}\right)^n\in k[x]
\]
is a solution to this differential equation.  Recall that $j\in k$ is
{\em supersingular} if and only if the corresponding elliptic curve 
over $k$ is supersingular. It follows that $j\neq 0, 1728$
is supersingular if and only if $u(j)=0$.

\section{The reduction of $X_0(p)$}

Let $X_0(p)$ denote the modular curve parameterizing (generalized)
elliptic curves with a $\Gamma_0(p)$-structure \cite{DelRap}. We may
and will identify $X_0(p)$ with the quotient curve $X(p)/G_0$, where
$G_0\subset G:={\rm PSL}_2(p)$ is the standard Borel group, modulo
$\pm I$. The $G$-cover $X(p)\to X(1)=\PP^1_j$ gives rise to a
non-Galois cover $X_0(p)\to\PP^1_j$ which we call the {\em
$j$-invariant}. We claim that the curve $X_0(p)$ has a {\em unique}
$\QQ$-model $X_0(p)_{\QQ}$ such that the $j$-invariant descends to a
cover $X_0(p)_{\QQ}\to\PP^1_{\QQ}$. (Of course, the existence of such
a $\QQ$-model is also a consequence of the modular interpretation of
$X_0(p)$.)

To prove the claim, let $G\hookrightarrow S_{p+1}$ denote the standard
permutation representation of $G$ coming from the natural action on
$\PP^1(\FF_p)$. Note that $G_0$ is precisely the stabilizer of
$\infty$. It is well known that the normalizer of $G$ in $S_{p+1}$ is
equal to ${\rm PGL}_2(p)=\Aut(G)$. In particular, the centralizer of
$G$ in $S_{p+1}$ is trivial. Let ${\bf C}=(3A,2A,pA)$ denote the
triple of conjugacy classes of $G$ corresponding to the $G$-cover
$X(p)\to X(1)$, as in the previous section. Let ${\bf C'}$ denote
the image of ${\bf C}$ in ${\rm PGL}_2(p)$. Then ${\bf C'}$ is
rigid and rational. Hence the claim follows from standard Galois
theory. 

For any field $K$ of characteristic zero we denote by
$X_0(p)_K$ the $K$-model of $X_0(p)$ obtained from
$X_0(p)_{\QQ}$ by base change. In \cite[VI.6.16]{DelRap} Deligne and
Rapoport prove the following result.

\begin{thm}
Let $X_0(p)_{\ZZ_p}$ denote the normalization of $\PP^1_{\ZZ_p}$ inside
$X_0(p)_{\QQ_p}$. Set $k:=\bar{\FF}_p$. 
\begin{enumerate}
\item
  The $\ZZ_p$-curve $X_0(p)_{\ZZ_p}$ is semistable (and stable if
  $g(X_0(p))>1$). 
\item
  The geometric special fiber
  $X_0(p)_k:=X_0(p)_{\ZZ_p}\otimes k$ is the union of two
  smooth curves $\Wb_0'$, $\Wb_0''$ of genus $0$. The induced map
  $\Wb_0'\to\PP^1_k$ (resp.\ $\Wb_0''\to\PP^1_k$)
  is an isomorphism (resp.\ purely inseparable of degree $p$).
  The components $\Wb_0'$ and $\Wb_0''$ intersect precisely in the {\em
  supersingular points}, i.e.\ the points of $X_0(p)_k$ with
  supersingular $j$-invariant. 
\item
  Let $x\in X_0(p)_k$ be a supersingular point, with $j$-invariant
  $j_x\in k$. If $j_x\equiv 0$ then $X_0(p)_{\ZZ_p}$ has a
  singularity of type $A_3$ at $x$. If $j_x\equiv 1728$ then $x$
  presents a singularity of type $A_2$. Otherwise, $x$ is a regular
  point of $X_0(p)_{\ZZ_p}$.
\end{enumerate}
\end{thm}

We give a proof of this theorem using the results of the previous
section.

\begin{Proof}
Let $R_0:=W(k)$ denote the ring of Witt vectors and $K_0$ the fraction
field of $R_0$. Note that
\[
    X_0(p)_{R_0} \;:=\; X_0(p)_{\ZZ_p}\otimes_{\ZZ_p}R_0
\]
is equal to the normalization of $\PP^1_{R_0}$ in the function field
of $X_0(p)_{K_0}$. Hence it suffices to prove the theorem with the
ring $\ZZ_p$ replaced by $R_0$. 

Let $K/K_0$ be the tame extension of degree $(p^2-1)/2$. By Corollary
\ref{xpcor} the ${\rm PSL}_2(p)$-cover $Y_K:=X(p)_K\to X_K:=\PP^1_K$
extends to a stable model $Y_R\to X_R$ over the ring  of integers $R$
of $K$. As before, we denote by $\Yb\to\Xb$ the special fiber of
$Y_R\to X_R$.

Set $W_R:=Y_R/G_0$. By \cite[Appendice]{Raynaudfest}, $W_R$ is a
semistable curve over $R$ with generic fiber $W_K=X_0(p)_K$. By
construction, we have a finite map $W_R\to X_R$. Let $\Wb$ denote the
special fiber of $W_R$. The formation of the quotient $Y_R/G_0$ does
not commute with base change, in general. However, the canonical map
$\Yb/G_0\to\Wb$ is a homeomorphism. Let $V$ be an irreducible
component of $\Yb$ and $U$ its image in $\Wb$. Let $D(V)\subset G$
(resp.\ $I(V)\subset G$) denote the decomposition group (resp.\ the
inertia group) of $V$. Then the natural map $V/(D(V)\cap G_0)\to U$ is
purely inseparable of degree $|I(V)\cap G_0|$. Using Corollary
\ref{xpcor} and the description of the stable reduction in Section 1,
we see that $\Wb$ has four different types of components (see Figure
\ref{x0ppic} for the case $p=29$).  In the description below we call
an irreducible component of $\Wb$ {\em horizontal} if it is mapped to
the central component $\Xb_0$ of $\Xb$ and {\em vertical} otherwies,
cf.\ Figure \ref{stablepic}.
\begin{itemize}
\item[(a)] A horizontal component $\Wb_0'$ which is the image of the
  (unique) horizontal component $\Yb_0\subset\Yb$ with decomposition
  group $G_0$. The natural map $\Wb_0'\to\Xb_0=\PP^1_k$ is an
  isomorphism.
\item[(b)] Another horizontal component $\Wb_0''$ which is the image of
  the set of horizontal components of $\Yb$ different from
  $\Yb_0$. Note that if $\Yb''\not=\Yb_0$ is such a horizontal
  component, with decomposition group $D''$, then $D''\cap G_0$ is a
  cyclic group of order $(p-1)/2$ (a split torus). Therefore, the
  natural map $\Wb_0''\to\Xb_0$ is purely inseparable of degree $p$.
\item[(c)] For $i\in \Bn\cup \Bp$ we let $\Wb_i$ be the image of the
  tail cover $\Yb_i$. Suppose that $i$ is supersingular, i.e.\ that
  the tail $\Xb_i$ intersects $\Xb_0=\PP^1_k$ in a supersingular
  point. Then $\Yb_i$ and --- hence $\Wb_i$ as well --- is
  connected. The curve $\Wb_i$ intersects each of the horizontal
  components $\Wb_0$ and $\Wb_0'$ in a unique point.
\item[(d)] With notation as in (c), suppose that $i$ is not
supersingular (in particular, $i\in \Bp$). Then $\Yb_i$ is not
connected, and the decomposition group of any connected component is
isomorphic to $\ZZ/p\rtimes\ZZ/2$ or $\ZZ/p\rtimes\ZZ/3$ (depending on
whether $x_i\equiv 1728$ or $x_i\equiv 0$). It follows that $\Wb_i$
has exactly two connected components $\Wb_i'$ and $\Wb_i''$.  The
curve $\Wb_i'$ (resp.\ $\Wb_i''$) intersects $\Wb_0'$ (resp.\
$\Wb_0''$) in a unique point. (See the end of the proof of Corollary
\ref{xpcor}.)
\end{itemize}

In the modular picture, the component $\Wb_0'$ (resp.\ $\Wb_0''$)
corresponds generically to ordinary elliptic curves together with a
subgroup scheme isomorphic to $\bmu_p$ (resp.\ to $\ZZ/p$), see
\cite[VI.6]{DelRap}. The vertical components have no easy modular
interpretation. Note that every component of $\Wb$ has genus $0$.

\begin{figure}[h] 
\begin{center} \unitlength3mm
  %%% Local Variables: 
%%% mode: latex
%%% TeX-master: "mcav"
%%% End: 

\begin{picture}(29,18)

\put(4,14){\line(1,0){22}}
\put(0,13){\makebox(4,2){$\Wb_0'$}}
\put(4,7){\line(1,0){22}}
\put(0,6){\makebox(4,2){$\Wb_0''$}}

\put(8,5){\line(0,1){4}}
\put(8,12){\line(0,1){4}}
\put(6,2){\makebox(4,2){$\scriptstyle j\equiv 1728$}}

\put(13,5){\line(0,1){11}}
\put(17,5){\line(0,1){11}}
\put(12,2){\makebox(6,2){$\scriptstyle j^2+2j-8\equiv 0$}}

\put(22,5){\line(0,1){11}}
\put(20,2){\makebox(4,2){$\scriptstyle j\equiv 0$}}

\end{picture}
  \caption{\label{x0ppic} The curve $\Wb$ for $p=29$}
\end{center}
\end{figure}

Let $W_R'$ be the curve obtained from $W_R$ by contracting all
vertical components of $\Wb$. Since these vertical components have
genus zero and self-intersection number $-1$ or $-2$, the curve $W_R'$
is semistable. In particular, $W_R'$ is a normal scheme. Recall that
the curve obtained from $X_R$ by contracting the tails $\Xb_i$ is
equal to $\PP^1_R$. Therefore the finite map $W_R\to X_R$ gives rise
to a finite map $W_R'\to\PP^1_R$. We conclude that $W_R'$ is the
normalization of $\PP^1_R$ inside the function field of
$W_K=X_0(p)_K$.

Let $\Gamma:=\Gal(K/K_0)\cong\ZZ/N$, and recall that $N=(p^2-1)/2$ is
prime to $p$. Since $W_K=X_0(p)_{K_0}\otimes_{K_0}K$, the group
$\Gamma$ acts naturally on $W_K$. Since $W'_R$ is the normalization of
$\PP^1_R$ in $W_K$ this action extends to the $R$-model $W_R'$ and
\begin{equation} \label{x0peq1}
    X_0(p)_{R_0} \;=\; W_R'/\Gamma.
\end{equation}
Since the order of $\Gamma$ is prime to $p$, formation of the quotient
commutes with base change. In particular, we have
\begin{equation} \label{x0peq2}
    X_0(p)_k \;=\; \Wb'/\Gamma.
\end{equation}
It is clear that the action of $\Gamma$ commutes with the map
$\Wb'\to\PP^1_k$ coming from the $j$-invariant. Since this map is an
isomorphism on one and purely inseparable of degree $p$ on the other
component, $\Gamma$  acts trivially on $\Wb'$. Now Parts (i) and (ii) of
the theorem follow from \eqref{x0peq2} and the description of the
components of $\Wb$ given above.

It remains to prove Part (iii) of the theorem. Let $x\in
X_0(p)_k=\Wb'$ be a supersingular point, with $j$-invariant $j_x\in
k$. Let $v$ denote the valuation on $K$, normalized such that
$v(p)=1$. Since $x$ is an ordinary double point of $\Wb'$, the local
ring of $X_0(p)_{R_0}$ at $x$ is of the form $R_0[[u,v\mid uv=\pi]]$,
with $\pi\in R_0$. We define the {\em thickness} of $X_0(p)_{R_0}$ at
$x$ as the rational number $e(x) \;:=\; v(\pi)$.  We have to show that
$e(x)$ is equal to $3$ if $j_x\equiv 1728$, equal to $2$ if $j_x\equiv
0$ and equal to $1$ otherwise.

We  assume that $j_x\not\equiv 0,1728$. The other cases follow
in the same manner. Let $x'\in\Wb'$ be the unique point lying above
$x$. Also, let $x''$ and $x'''$ be the two double points on $\Wb$
lying above $x'$. We assume that $x''\in\Wb_0'$ and
$x'''\in\Wb_0''$. It is easy to see that
\begin{equation} \label{x0peq3}
          e(x) \;=\; e(x') \;=\; e(x'')+e(x''').
\end{equation} 
Let $y'',y'''\in\Yb$ be points above $x'',x'''$. Recall that
$W_R=Y_R/G_0$. By \cite[Appendice]{Raynaudfest}, we have
$e(x'')=|{\rm Stab}_{G_0}(y'')|\cdot e(y'')$ and 
$e(x''')=|{\rm Stab}_{G_0}(y''')|\cdot e(y''')$. Using \eqref{x0peq3}
and our knowledge of the $G$-action on $\Yb$, we see that
\begin{equation} \label{x0peq4}
     e(x) \;=\; \frac{p(p-1)}{2}\cdot e(y'')+\frac{p-1}{2}\cdot e(y''').
\end{equation}
It is proved in \cite{bad} that the thickness of any double point $y\in
\Yb$ is equal to $(h_i(p-1))^{-1}$ if $y$ lies on the tail $\Yb_i$ and
$h_i$ denotes the conductor of the tail cover $\Yb_i\to\Xb_i$ at
$\infty$. Since we assumed that $j_x\not\equiv 0,1728$ we have
$y'',y'''\in\Yb_i$ for some $i\in\Bn$. By Corollary \ref{xpcor} we have
$h_i=(p+1)/2$ for $i\in\Bn$. With \eqref{x0peq4} we conclude that
\[
   e(x) \;=\; \frac{p}{p+1} + \frac{1}{p+1} \;=\; 1.
\]
This finishes the proof of the theorem.
\end{Proof}
 
One can analyze the natural $R_0$-model of $X_1(p)$ in an entirely
analogous manner. This would give a new proof of the result of
\cite{DelRap} that $J_1(p)/J_0(p)$ has good reduction over
$\QQ(\zeta_p)$. One can also derive this result directly from our
knowledge of the stable reduction of $X(p)$.

%%% Local Variables: 
%%% mode: latex
%%% TeX-master: "mcav"
%%% End: 

\section{Three point covers with Galois group $\SL_2(p)$}\label{legosec}

In this section we compute the stable reduction of all
$\SL_2(p)$-covers of $\PP^1$ branched at three points. In particular,
we obtain an explicit formula for the number of such covers with good
reduction (Theorem \ref{reductionthm}). To prove this we use the
results of the previous sections. In addition, we use the fact that
some ${\rm SL}_2(p)$-covers have a modular interpretation as a {\em
  Hurwitz space}.

Let us give an outline of the proof of Theorem \ref{reductionthm}.  We
first show that the $\SL_2(p)$-covers of $\PP^1$ branched at three
points of order $p,p,n$ are essentially Hurwitz spaces. These covers
have bad reduction to characteristic $p$; their stable reduction is
described in \cite{IreneAux}. By considering the (unique) primitive
tail of the stable reduction, we obtain one primitive tail cover for
every noncentral conjugacy class of $\SL_2(p)$ whose order is prime to
$p$.  Proposition \ref{sigmaprimprop} states that every primitive
$\SL_2(p)$-tail cover with group ${\rm SL}_2(p)$ arises in this
way. Then we show (Proposition \ref{sigmanewprop}) that the only new
tail cover that occurs is the cover constructed in Lemma
\ref{AL}. This implies that the special deformation datum of an
$\SL_2(p)$-cover branched at three points is hypergeometric
(Definition \ref{hgdef}). Therefore, Proposition \ref{liftprop} allows
to construct all $\SL_2(p)$-covers with bad reduction at $p$, by
lifting.  This also gives a formula for the number of covers with good
reduction. Passing to the quotient modulo $-I$, one obtains similar
result for three point covers with Galois group $\PSL_2(p)$.

In this section, $R$ is a complete mixed characteristic discrete
valuation ring which contains $\sqrt{p^*}$ and whose residue field is
an algebraically closed field $k$ of characteristic $p$. We let $K$
denote the fraction field of $R$.

\begin{lem}\label{Hurwitzspacelem} 
  Let $D\in\{C(l),\tilde{C}(l)\}$. Let $f:Y\to\PP^1_K$ be an
  $\SL_2(p)$-Galois cover with class vector ${\bf C}_1=(pA, pA, D)$ or
  ${\bf C}_2=(pA, pB, D)$. Let $\Yb\to\Xb$ be the stable reduction of
  $f$. The curve $\Xb$ has exactly one primitive
  tail, with ramification invariant
  \[\renewcommand{\arraystretch}{2.0}
    \sigma \;=\;\; \left\{\begin{array}{cc}
      \displaystyle{\frac{p-1-2l}{p-1}} &\quad\text{if\ $ D=C(l)$},\\
      \displaystyle{\frac{p+1-2l}{p-1}} &\quad\text{if\ $D=\tilde{C}(l)$}.
    \end{array}\right.
  \]
  The inertia group of the primitive tail has order
  $p(p-1)/\gcd(p-1,l)$ if $ D=C(l)$ and $p(p-1)/\gcd(p-1, l-1)$
  otherwise.  The new tails have ramification invariant
  $\sigma=(p+1)/(p-1)$ and an inertia group of order $p(p-1)$.
\end{lem}

For $\SL_2(p)$-covers with ramification of order
$n$, one can prove a strengthening of Proposition
\ref{linrigprop}, cf.\ \cite[Lemma 3.29]{Voelklein}. Either there
exists a cover with class vector ${\bf C}_1$ or there exists a cover
with class vector ${\bf C}_2$. In both cases, there is a unique such
cover. The outer automorphism of $\SL_2(p)$ conjugates the cover with
class vector ${\bf C}_1$ (resp.\ ${\bf C}_2$) to a cover with class
vector $(pB, pB, D)$ (resp.\ $(pB, pA, D)$). Therefore, in both cases,
there are exactly two covers with ramification $(p,p, D)$, up to isomorphism.

\begin{Proof} The case $D=C(l)$ is treated in \cite[Section
  8]{IreneAux}. The idea of the proof is as follows. Define
  $m=(p-1)/\gcd(p-1,l)$. For $\lambda\in\PP^1-\{0,1,\infty\}$, let
  $W_\lambda$ be a connected component of the nonsingular projective curve
  corresponding to 
  \[
       w^{p-1} \;=\; x^l(x-1)^l(x-\lambda)^{p-1-l}.
  \]
  This defines a family of $m$-cyclic covers
  $g_\lambda:W_\lambda\to\PP^1$, branched at four points. Fix an
  injective character $\chi:\ZZ/m\to\FF_p^\times$ and let
  $N=\ZZ/p\rtimes_\chi \ZZ/m$.  Define the Hurwitz space $\HH$
  parameterizing $N$-Galois covers branched at $0,1,\lambda, \infty$
  of order prime to $p$ which factor via $g_\lambda$. Define
  $\HH\to\PP^1_\lambda$ by mapping an $N$-cover to the branch point
  $\lambda$ and let $f:\tilde{\HH}\to\PP^1_\lambda$ its Galois
  closure. The Galois group of $f$ is either $\SL_2(p)$ or
  $\PSL_2(p)$. If it is $\SL_2(p)$, the class vector of $f$ is $(D_1,
  D_2 , D)$, with $D_1, D_2\in\{pA, pB\}$. Otherwise, there exists a
  unique lift of $f$ to a cover with such a class vector. Therefore,
  $f$ is isomorphic to the cover of the statement of the proposition.
  In \cite{IreneAux}, the stable reduction of $f$ is computed, by
  using the moduli interpretation of $\HH$.

  The proof for $D=\tilde{C}(l)$ is analogous. (Take $m|(p+1)$ and
  $N\simeq (\ZZ/p)^2\simeq \ZZ/m$.)
\end{Proof}

The following lemma is due to R.\ Pries (unpublished).

\begin{lem}\label{liftlem} 
  Let $G$ be a finite group. Let $\varphi:W\to \PP^1_k$ be a new
  $G$-tail cover with ramification invariant $\sigma$ and inertia
  group $I$ of order $pn$. Suppose $1<\sigma <2$. Then there exists a
  $G$-cover $f:Y\to\PP^1_K$ with signature $(p,p,n)$ such that
  $\varphi$ occurs as the unique new tail of the stable reduction of
  $f$.
\end{lem}

If $\sigma=2$ then $\varphi$ can be lifted to a cover in
characteristic zero which is branched at three points of order $p$,
\cite[Section 3.2]{RRR}.

\begin{Proof} Write $n=m_1m_2$, where $m_2$ is the
order of the prime-to-$p$ part of the center of $I$. Then $m_1$
divides $p-1$ and $\sigma-1=a/m_1$ for some $0<a<m_1$. Write
$\alpha=(p-1)/m_1$. Choose $\lambda\in k-\{0,1\}$ such that
$(1-\lambda)^{\alpha a}=1$. Note that $\sigma\neq p/(p-1)$, since
conductors are prime to $p$. Therefore, $\alpha a\neq 1$ and it is
possible to choose $\lambda$ as asserted.
  
Let $g_0:\Zb_0\to \Xb_0$ be the $m_1$-cyclic cover where $\Zb_0$ is
defined by ${z}^{m_1}=x^{-a}(x-\lambda)^{a}$. Let $g:\Zb\to\PP^1_k$ be the
$n$-cyclic cover branched only at $0$ and $\lambda$ which has $g_0$ as
a quotient. Define the differential
\[
   \omega \;=\; \frac{\epsilon\,z\,{\rm d}x}{x(x-1)}
\]
on $\Zb$, with $\epsilon$ as in \eqref{epsiloneq}. A straightforward
computation, as in the proof of Proposition \ref{hgprop}, shows that
$(g,\omega)$ is a special deformation datum. There is one primitive
critical point at zero and one new critical point at $\lambda$.

Define a primitive $I$-tail cover $\psi':\Yb'\to\PP^1_k$ which is
totally ramified at $\infty$ with ramification invariant
$\sigma=(m_1-a)/m_1$ and branched of order $n$ at 0. Let
$\psi:=\Ind_I^G\psi'$. One checks that the datum
$(g,\omega,\psi,\varphi)$ satisfies the conditions of Proposition
\ref{liftprop}. We conclude that there exists a $G$-Galois cover
$f:Y\to\PP^1_K$ which is branched at three points with ramification of
order $(p,p,n)$ such that its stable reduction gives rise to
$(g,\omega,\psi,\varphi)$. This proves the lemma.
\end{Proof}

It follows from rigidity that the special deformation datum associated
to an ${\rm SL}_2(p)$-cover with bad reduction is uniquely determined
by the invariants $h_i$ and $m_i$. On the other hand, the special
deformation datum $(g,\omega)$ in the proof of Lemma \ref{liftlem}
depends also on the choice of $\lambda\in k-\{0,1\}$ with
$(1-\lambda)^{\alpha a}=1$. It follows that we have $\alpha a=2$ in
case $G=\SL_2(p)$.  This statement is also a consequence of the next
proposition.

\begin{prop}\label{sigmanewprop}
  Suppose that $G$ is either $\SL_2(p)$ or $\PSL_2(p)$.  Let
  $\varphi:W\to \PP^1_k$ be a new $G$-tail cover with ramification
  invariant $\sigma$. Suppose that $1<\sigma\leq 2$. Then
  $\sigma=(p+1)/(p-1)$. The order of the inertia group is $p(p-1)$ if
  $G=\SL_2(p)$ and $p(p-1)/2$ if $G=\PSL_2(p)$. Furthermore, the cover
  $\varphi$ is unique, up to isomorphism.
\end{prop}

\begin{Proof} 
  If there exists a new $\SL_2(p)$-tail cover $\varphi$ with $\sigma\neq
  2$, Lemma \ref{liftlem} implies that $\varphi$ occurs as the
  restriction to a new tail of the stable reduction of a ${\rm
    SL}_2(p)$-cover with signature $(p,p,n)$ for some $n$
  prime-to-$p$.  Therefore Lemma \ref{Hurwitzspacelem} implies that
  $\sigma=(p+1)/(p-1)$.
  
  Similarly, one shows that every $\SL_2(p)$-cover of $\PP^1_k$
  branched at one point, with $\sigma=2$, can be lifted to a cover of
  $\PP^1_K$ which is branched at three points of order $p$. The
  proposition follows in this case from Proposition \ref{redmodprop}
  and Proposition \ref{linrigprop}, cf.\ \cite[Corollary 3.2.2]{RRR}.
  
  Since every $\PSL_2(p)$-cover of $\PP^1_K$ with ramification of
  order $p,p,n'$ can be  lifted to an $\SL_2(p)$-cover with
  ramification of order $p,p,n$ (where $n$ is $n'$ or $2n'$), the
  proposition also follows for $G=\PSL_2(p)$.
 \end{Proof}

\begin{cor}
  Let $G$ be $\SL_2(p)$ or $\PSL_2(p)$.  Let $f:Y\to\PP^1_K$ be a
  $G$-Galois cover branched at three points. Suppose that $f$ has bad
  reduction. Then the special deformation datum corresponding to $f$
  is hypergeometric.
\end{cor}

\begin{Proof}$f:Y\to\PP^1_K$ be a
  $G$-Galois cover branched at three points. Suppose that $f$ has bad
  reduction.  Proposition \ref{sigmanewprop} implies that
  $\sigma_i=(p+1)/(p-1)$ for $i\in\Bn$. Therefore the corollary
  follows from Definition \ref{hgdef}.
\end{Proof}

\begin{prop}\label{sigmaprimprop}
\begin{itemize}
\item[(i)] Let $\varphi:W\to \PP^1_k$ be a (possibly nonconnected)
  primitive $\SL_2(p)$-tail cover defined over $k$, with ramification
  invariant $\sigma$.  Suppose that $0< \sigma<1$. Denote by $D$ the
  conjugacy class in $\SL_2(p)$ of the canonical generator of a point
  of $W$ above 0, with respect to the fixed roots of unity.  Then
  \[\renewcommand{\arraystretch}{2.0}
    \sigma \;=\quad \left\{\begin{array}{cc}
    \displaystyle{\frac{p-1-2l}{p-1}} &\quad\text{if\ $D=C(l)$},\\
    \displaystyle{\frac{p+1-2l}{p-1}} &\quad\text{if\ $D=\tilde{C}(l)$}.
    \end{array}\right.
  \]
\item[(ii)] The decomposition group of an irreducible component of $W$
is $\SL_2(p)$ if $D=\tilde{C}(l)$ and is a semidirect product of order
$p(p-1)/gcd(p-1,l)$  if $C=C(l)$.
\end{itemize}
\end{prop}

\begin{Proof}  
  Write $\sigma=a/(p-1)$.  The normalizer $N_G(P)$ of a Sylow
  $p$-subgroup $P$ in $\SL_2(p)$ has order $p(p-1)$. Since the center
  of $N_G(P)$ has order two, it follows that $a $ is even.  Therefore
  $(p-1-a )/2$ is a positive integer.
 
  Proposition \ref{hgprop} (ii) shows that there exists a
  hypergeometric deformation datum $(g',\omega')$ of signature $(0,
  0,a)$, where $g':Z'_k\to\PP^1_k$ is a $(p-1)/2$-cyclic cover,
  given by
  \[ 
     {z'}^{(p-1)/2} \;=\; \prod_{i\in\Bn} (x-\tau_i).
  \]
  We define $g:Z_k\to\PP^1_k$ to be the $(p-1)$-cyclic
  cover given by
  \[ 
     z^{p-1} \;=\; \prod_{i\in\Bn} (x-\tau_i).
  \]
  Let $\omega$ be the pull back of $\omega'$ to $Z_k$.  Then
  $(g,\omega)$ is a special deformation datum.  Note that $|\Bp|=1$
  and $|\Bn|=(p-1-a)/2$, with $a=p-1-2l$ if $D=C(l)$ and $a= p+1-2l$
  if $D=\tilde{C}(l)$. The primitive critical point is at $\infty$,
  where the cover $g$ is branched of order $(p-1)/(p-1,l)$ if $D=C(l)$
  and of order $(p-1)/(p-1,l-1)$ otherwise.
  
  Let $h$ be the $\SL_2(p)$-tail cover constructed in Lemma \ref{AL}.
  For every $i\in\Bn$, let $\fb_i$ be a copy of $h$.  Proposition
  \ref{liftprop} shows that there exists an $\SL_2(p)$-cover $f:Y
  \to\PP^1_K$ whose stable reduction gives rise to the special
  deformation datum $(g , \omega )$, the new tail covers $\fb_i$ and
  the primitive tail cover $\varphi$. By construction, the class
  vector of $f$ is $(D_1, D_2, D)$ with $D_1, D_2\in\{pA, pB\}$.
  Therefore $f$ is isomorphic to a cover described in Lemma
  \ref{Hurwitzspacelem}. Part (i) of the proposition follows from
  this.
  
  Let $N$ denote the decomposition group of an irreducible component
  of $W$.  Suppose $D =\tilde{C}(l)$, for some $l<(p-1)/2$. Then $N$
  contains an element of order $p$ and an element of order
  $(p+1)/\gcd(p+1,l)>2$.  Therefore, $N=\SL_2(p)$, \cite[Section
  II.8]{Huppert}.
  
  Next suppose that $D=C(l)$. Put $n_l=(p-1)/\gcd(p-1,l)$. It is easy
  to construct a primitive tail cover with Galois group $\ZZ/p\rtimes
  \ZZ/n_l$ which is branched at zero of order $n_l$ and has
  ramification invariant $\sigma=(p-1-2l)/(p-1)$. Let $\varphi'$ be
  the primitive $\SL_2(p)$-tail cover obtained by induction. Using
  again Proposition \ref{liftprop}, we see that there exists an
  $\SL_2(p)$-cover $f':Y' \to\PP^1_K$ whose stable reduction gives
  rise to the special deformation datum $(g,\omega )$, the new tail
  covers $\fb_i$ and the primitive tail cover $\varphi'$. Note that
  $f'$ has the same branch cycle description as $f$. It follows from
  Proposition \ref{linrigprop} that $f$ and $f'$ are outer isomorphic,
  i.e.\ become isomorphic as ${\rm SL}_2(p)$-covers after twisting by
  an outer automorphism of ${\rm SL}_2(p)$. We conclude that $N$ is an
  extension of $\ZZ/n_l$ by $\ZZ/p$.
\end{Proof}

\begin{thm}\label{reductionthm}
  Let ${\bf C}=(C_1, C_2, C_3)$ be a triple of conjugacy classes of
  $\SL_2(p)$, where we suppose that the elements of the $C_i$ are not
  contained in the center.
\begin{itemize}
\item[(a)] If $C_i$ contains an element of order $p$ for some $i$,
   all covers with class vector ${\bf C}$ have bad reduction.
 \item[(b)] Suppose $C_i\in\{C(l), \tilde{C}(l)\}$ and put
   $a_i=p-1-2l$ if $C_i=C(l)$ and $a_i= p+1-2l$ otherwise. 
   Write 
   \[
       \Ni^\bad_3({\bf C}) \;=\; \{\, f\in\Niin_2({\bf C})\, |\, 
           \mbox{ $f$ has bad reduction.}\,\}.
   \]
  Then
  \[
       |\Ni^\bad_3({\bf C})| \;=\;
        \begin{cases}
          \;2 & \text{ if }a_1+a_2+a_3< p-1,\\
          \;2 & \text{ if }a_1+a_2+a_3=p-1 \text{ and }C_i=
           \tilde{C}(l)\text{ for some }i,\\
          \;0 &\text{ otherwise.}
        \end{cases}.
  \]
  \end{itemize}
\end{thm}

\begin{Proof} Let $f:Y\to\PP^1_K$ be an $\SL_2(p)$-cover with class 
  vector ${\bf C}=(C_1, C_2, C_3)$.  If $p$ divides one of the
  ramification indices of $f$, then $f$ has bad reduction, so (a) is
  obvious.
  
  Suppose that $p$ does not divide the ramification indices of $f$. If
  $f$ has bad reduction, then (\ref{vceq}) implies that
  $a_1+a_2+a_3\leq p-1$. If $a_1+a_2+a_3=p-1$, there are no new tails.
  Since $\Yb$ is connected, at least one of the primitive tails has
  decomposition group $\SL_2(p)$.  Therefore Proposition
  \ref{sigmaprimprop} implies that $C_i=\tilde{C}(l)$, for some $i$.
  This shows that if there exists a cover with bad reduction, then
  either $a_1+a_2+a_3<p-1$ or $a_1+a_2+a_3=p-1$ and $C_i=\tilde{C}(l)$
  for some $i$.
  
  Now suppose that $C_1, C_2, C_3$ are conjugacy classes of noncentral
  elements of order prime to $p$ such that the condition of the
  theorem is satisfied. Let $a_1, a_2, a_3$ be as in the statement of
  the theorem. Recall that $a_i$ is even and different from $ p-1$. Let
  $(g',\omega')$ be the hypergeometric deformation datum of signature
  $(a_1, a_2, a_3)$ constructed in Proposition \ref{hgprop} (ii). The
  cover $g':Z'_k\to\PP^1_k$ is the $(p-1)/2$-cyclic cover defined by
  \[ 
     {z'}^{(p-1)/2} \;=\; x^{a_1/2}(x-1)^{a_2/2}
        \prod_{i\in\Bn}(x-\tau_i), \quad (z,x)\mapsto x.
  \]
  We define a $(p-1)$-cyclic cover $g:Z_k\to\PP^1_k$,
  factoring through $g'$ by
  \[
     z^{p-1} \;=\; x^{b_1}(x-1)^{b_2}
        \prod_{i\in \Bn}(x-\tau_i),\quad (z,x)\mapsto x,
  \]
  where $b_i=a_i/2+(p-1)/2$, i.e.\ $b_i=p-1-l$ if $C_i=C(l)$ and
  $b_i=p-l$ otherwise. Then $(g,\omega)$ is a special deformation
  datum.
  
  We define a set of tail covers, as follows. For $i\in\Bn$, we define
  $\fb_i\simeq h$, where $h$ is the new $\SL_2(p)$-tail cover of
  $\PP^1_k$ with $\sigma=(p+1)/(p-1)$ from Lemma \ref{AL}.  For
  $i\in\Bp=\{1,2,3\}$, we let $\fb_i$ be the primitive $\SL_2(p)$-tail
  cover with $\sigma=a_i/(p-1)$ whose existence is guaranteed by Lemma
  \ref{Hurwitzspacelem}. As in the proof of Lemma \ref{liftlem}, one
  checks that the datum $(g,\omega;\fb_i)$ satisfies the conditions of
  Proposition \ref{liftprop}. Hence there exists an $\SL_2(p)$-cover
  $f:Y\to\PP^1_K$ with class vector ${\bf C}=(C_1, C_2, C_3)$ whose
  stable reduction gives rise to the datum $(g,\omega;\fb_i)$.
  
  Since the outer automorphism group of $\SL_2(p)$ has order two, the
  number of liftings, up to isomorphism, is at least two.
  Proposition \ref{linrigprop} implies that the number of elements of
  $\Niin_3({\bf C})$ is at most two. We conclude that the number of
  liftings is exactly two. This proves that all $\SL_2(p)$-covers over $K$ with
  class vector ${\bf C}$ have bad reduction to characteristic $p$.
\end{Proof}

\vspace{4ex}
\noindent{\small
  Institut f\"ur Experimentelle Mathematik \hfill
    Max--Planck--Institut f\"ur Mathematik\\
  Universit\"at GH Essen\hfill
    Vivatsgasse 7\\
  Ellernstr. 29\hfill
    53111 Bonn, Germany\\
  45326 Essen, Germany \hfill
    wewers@mpim-bonn.mpg.de\\
  bouw@exp-math.uni-essen.de}
\vspace{3ex}

\end{document}